\documentclass[11pt]{amsart}
\usepackage[margin=1.08in]{geometry}
\usepackage[T1]{fontenc}
\usepackage{amssymb,microtype}
\usepackage{tikz}
\usepackage[hidelinks]{hyperref}
\newcommand{\manuscriptkeywords}{Garnett--Jones distance theorem, rectangular BMO,
  little bmo, exponential integrability, coordinate maximal operators,
  two-weight estimates}
\hypersetup{
  pdftitle={The Garnett-Jones distance theorem for rectangles},
  pdfauthor={Morten Nielsen and Hrvoje \v{S}iki\'{c}},
  pdfkeywords={\manuscriptkeywords}
}
\newcommand{\one}{\mathbf 1}
\newcommand{\av}[2]{\langle #1\rangle_{#2}}
\newcommand{\pair}[2]{\langle #1,#2\rangle}
\newcommand{\BMO}{\mathrm{BMO}_{2,\mathcal R}}
\newcommand{\D}{\mathcal D}
\newcommand{\rev}{\mathrm{rev}}
\newtheorem{theorem}{Theorem}[section]
\newtheorem{lemma}[theorem]{Lemma}
\newtheorem{proposition}[theorem]{Proposition}
\allowdisplaybreaks[2]
\title{The Garnett--Jones distance theorem for rectangles}
\author[M. Nielsen]{Morten Nielsen}
\address{Department of Mathematical Sciences, Aalborg University,
  DK-9220 Aalborg {\O}st, Denmark}
\email{\href{mailto:mnielsen@math.aau.dk}{mnielsen@math.aau.dk}}
\author[H. \v{S}iki\'{c}]{Hrvoje \v{S}iki\'{c}}
\address{Department of Mathematics, Faculty of Science, University of Zagreb,
  HR-10000 Zagreb, Croatia}
\email{\href{mailto:hsikic@math.hr}{hsikic@math.hr}}
\begin{document}
\begin{abstract}
For the family of all axis-parallel rectangles in $\mathbb R^d$, $d>1$,
we characterize the distance of a real- or complex-valued function of
rectangular bounded mean oscillation (little bmo) to $L^\infty$ by the
reciprocal of its critical exponential integrability exponent. For each
$\alpha>0$ at which the rectangular exponential moment of $f$ is finite,
we construct a bounded correction so that the residual has mean absolute
oscillation bounded by a dimension-dependent constant times
$1/\alpha$ on every rectangle. The proof
bounds products of coordinate maximal operators between two sets in
terms of the supremum of their density product over common test
rectangles. Selecting a longest interval along each path lets us remove
repeated visits to a coordinate, so the density exponent is independent
of the number of factors. The maximal-operator estimate holds with every
positive density exponent below $1/6$ in two dimensions and below
$1/24$ in three dimensions.
\end{abstract}
\maketitle

\section{Introduction and main results}
Garnett and Jones characterized the distance from a BMO function to the
bounded functions in terms of uniform exponential integrability of its
oscillation on cubes \cite{gj}. Varopoulos gave a proof based on Brownian
motion \cite{varopoulos}. The original proof also led to a geometric
separation problem. If each cube contains a sufficiently small proportion
of at least one set from a given finite family, Uchiyama constructs a
partition of unity whose members have small BMO norms and vanish almost
everywhere on the corresponding sets \cite{uchiyama}.

A weight-based proof on cubes uses Jones's factorization of an $A_p$
weight as $w=uv^{1-p}$ with $u,v\in A_1$ and $1<p<\infty$
\cite{jones-factorization}. Coifman and Rochberg showed that each $A_1$
factor can be written as $(Mg)^\delta$, with $M$ the Hardy--Littlewood
maximal operator and $0<\delta<1$, times a function bounded above and
away from zero. Their maximal-power estimate gives an $A_1$ constant for
$(Mg)^\delta$ depending only on $\delta$ and the dimension
\cite{coifman-rochberg}. Applying the estimate with exponent $1/2$ and
using Jensen's inequality gives a dimensional BMO bound for $\log Mg$,
and hence for $\delta\log Mg$. Taking logarithms in the factorization
then splits $\log w$ into a function of BMO norm $O_{d,p}(1)$ and a
bounded function. For rectangles the factorization step survives: Rubio
de Francia's maximal-operator proof of factorization applies to the
rectangle basis \cite{rubio}. Soria
showed that powers of the strong maximal function can fail the
rectangular $A_1$ condition \cite{soria}. Nielsen and \v{S}iki\'c
related the upper distance estimate for a class of bases to an
$A_2$-decomposition property and left the case of all axis-parallel
rectangles open \cite{ns}. We prove the
characterization for rectangles in each fixed dimension. Instead of
factorizing, we build the bounded correction directly from a two-weight
estimate for products of coordinate maximal operators.

The single-average oscillation used here is the Euclidean counterpart of
little BMO on the polydisk. On the distinguished boundary of the bidisk,
Ferguson and Sadosky characterized little BMO by $L^2$ boundedness of the
commutators with the two coordinate Hilbert transforms
\cite{ferguson-sadosky}. Little BMO is a proper subspace of the product
BMO space of Chang and Fefferman, the dual of real product $H^1$
\cite{chang-fefferman}, which is not considered here.

Let $\mathcal R$ consist of all products of $d$ bounded intervals of
positive length in $\mathbb R^d$. Write $\av gR=|R|^{-1}\int_R g$
and $f_R=\av fR$. For $p=1,2$ set
\[
 \|f\|_{\mathrm{BMO}_{p,\mathcal R}}
    =\sup_{R\in\mathcal R}\bigl(\av{|f-f_R|^p}R\bigr)^{1/p},
 \qquad
 K_\alpha(f)=\sup_{R\in\mathcal R}\av{e^{\alpha|f-f_R|}}R.
\]
The BMO seminorm and the exponential moment both use the average $f_R$
over the whole rectangle. The seminorms $\|\cdot\|_{\mathrm{BMO}_{1,\mathcal R}}$ and
$\|\cdot\|_{\mathrm{BMO}_{2,\mathcal R}}$ are equivalent with constants
depending only on $d$. Cauchy--Schwarz bounds the first by the second.
For the reverse bound, use $t^2\le 2e^t$ for $t\ge0$ in the rectangular
John--Nirenberg estimate in Section~\ref{converse}, which also treats
complex-valued functions.

By Jensen's inequality $K_\alpha(f)\ge1$, and
$K_\alpha(f)$ increases with $\alpha$. Call $\alpha>0$ \emph{admissible}
for $f$ if $K_\alpha(f)<\infty$; admissibility is downward closed. Write
$a_*(f)=\sup\{\alpha>0:K_\alpha(f)<\infty\}$ and
$\operatorname{dist}_1(f,L^\infty)=\inf_{h\in L^\infty}
\|f-h\|_{\mathrm{BMO}_{1,\mathcal R}}$.
For real $f$ the infimum may be taken over real $h$, since taking the
real part can only decrease the oscillation seminorm of $f-h$.

\begin{theorem}\label{distance}
Fix $d>1$. There are constants $c_d,C_d>0$ depending only on $d$ with
the following properties. If $f\in L^1_{\rm loc}(\mathbb R^d)$ is real
and $\alpha$ is admissible for $f$, there exists a real
$h\in L^\infty(\mathbb R^d)$ such that
\[
 \|h\|_\infty\le \alpha^{-1}\log K_\alpha(f),\qquad
 \|f-h\|_{\mathrm{BMO}_{1,\mathcal R}}\le C_d/\alpha .
\]
For complex-valued $f\in L^1_{\rm loc}(\mathbb R^d)$ and admissible
$\alpha$, there exists a bounded complex-valued $h$ such that
\[
 \|h\|_\infty\le \frac{\sqrt{2}}{\alpha}\log K_\alpha(f),\qquad
 \|f-h\|_{\mathrm{BMO}_{1,\mathcal R}}\le \frac{\sqrt{2}C_d}{\alpha}.
\]
For real- or complex-valued $f\in\mathrm{BMO}_{1,\mathcal R}$
one has $a_*(f)>0$ and
\[
 \frac{c_d}{a_*(f)}\le \operatorname{dist}_1(f,L^\infty)
   \le \frac{\sqrt2\,C_d}{a_*(f)},
 \qquad 1/\infty=0 .
\]
In particular the distance is zero exactly when every $\alpha>0$ is
admissible.
\end{theorem}

To obtain the upper distance bound, we construct a bounded correction on
each finite grid and pass to a single correction on all rectangles.
The grid construction requires bounds for products of coordinate
maximal operators tested against indicators of two sets.

We work in three models. On $\mathbb R^d$ we use Lebesgue measure and
ordinary intervals. On a product of finite lines we use counting measure.
On an \emph{ordinary grid} the intervals of an axis are the sets of
consecutive cells. On a \emph{dyadic grid} each axis has $2^m$ cells for
some $m$, and its intervals are the nodes of a fixed binary interval tree
whose root is the entire axis, whose leaves are singletons, and whose two
children at each split have equal size. Every parallel fiber uses the
same tree. A \emph{rectangle} in a model is a product of intervals of
that model.

Let $M_i h$ be the supremum of the averages of $|h|$ over intervals of
the model through the point along axis $i$, with the other coordinates
held fixed. For a word $w=i_1\cdots i_k$ put $M_w=M_{i_1}\cdots M_{i_k}$;
operators act on functions from right to left, and the empty word is the
identity. For a pair of sets $E,F$ write
\[
 \theta=\theta(E,F)=\sup_R\av{\one_E}R\av{\one_F}R,\qquad
 S=\sqrt{|E||F|},
\]
where the supremum runs over the rectangles of the model under
consideration. Thus $\theta$ measures how dense the two sets can be on a
common test rectangle.

\begin{theorem}\label{word}
For $d\ge2$ define $g_2=1/6$ and $g_{d+1}=(6+3/g_d)^{-1}$, so that
$g_3=1/24$. For each $d>1$ and $0<\eta_d<g_d$ there is $A_d<\infty$
such that
\[
 \pair{\one_E}{M_w\one_F}\le A_d^k\,\theta(E,F)^{\eta_d}\sqrt{|E||F|}
\]
for every word $w$ of length $k\ge1$ and all sets $E,F$ of finite measure,
where $\pair gh$ denotes the integral or sum of $gh$.
The bound holds on dyadic grids, ordinary grids, and $\mathbb R^d$.
The constant $A_d$ depends only on $d$ and $\eta_d$; it is uniform in
the word, the sets, and the grid sizes.
\end{theorem}

The exponents $g_d$ are not optimized. The proof of
Theorem~\ref{word} starts with words that visit each coordinate once.
Intersecting rectangles can be longer in different coordinates, as in
C\'ordoba and Fefferman's covering proof of the strong maximal theorem
\cite{cordoba-fefferman}. In the plane, a crossing of horizontal and
vertical level sets yields a single rectangle with a quantitative lower
bound on the product of the two densities. An induction adds one
coordinate at a time (Sections~\ref{planesection}--\ref{dimension}). For
a repeated coordinate, we select a longest interval along each averaging
path. Two applications of H\"older's inequality let us keep this step
and remove the other visits to that coordinate (Section~\ref{compress}).
Random shifts transfer the dyadic estimate to ordinary grids, and strong
$L^2$ convergence gives the continuum word estimate
(Section~\ref{ordinary}).

For Theorem~\ref{distance}, we use the ordinary finite-grid estimate to
bound weighted products of averaging operators (Section~\ref{weights}).
For real $f$ and admissible $\alpha$, apply these bounds to the weights
$e^{\pm\alpha f_n}/K_\alpha(f)$, where $f_n$ is the cell-average
function. This gives a grid function of uniformly bounded mean
oscillation within $\log K_\alpha(f)$ of $\alpha f_n$.
The construction uses a positive iteration of Rubio de Francia type
\cite{rubio} on pairs of grid functions (Section~\ref{barrier}). The
corresponding corrections have a weak-$*$ limit that works on all
rectangles.
The lower distance bound follows by applying interval exponential
integrability across the coordinates (Section~\ref{converse}).

\section{Preliminaries}\label{prelim}
\paragraph{Interval estimates.}
Let $M$ be the one-dimensional maximal operator of a model. On a dyadic
line the maximal intervals on which $|h|$ has average greater than $t$
are disjoint, so $|\{Mh>t\}|\le\|h\|_1/t$. For ordinary intervals on a
finite line, select disjoint witnessing intervals in decreasing order of
length. Every discarded interval meets a selected one of at least its
length and lies in its concentric triple. Hence the triples of the
selected intervals cover the level set, and
$|\{Mh>t\}|\le3\|h\|_1/t$. On the real line the finite covering
argument applies to compact subsets of the level set by enlarging
witness intervals slightly, and inner regularity gives the full
estimate. For $1<p<\infty$, the bound
$\{Mh>t\}\subseteq\{M(|h|\one_{\{|h|>t/2\}})>t/2\}$,
integrated against $p t^{p-1}\,dt$, yields
$\|Mh\|_p\le C_p\|h\|_p$. Summing or integrating along fibers gives
the coordinate $L^p$ bounds. Each coordinate maximal operator is
sublinear, is monotone on nonnegative inputs, and preserves nonnegative
constants. On nonnegative inputs it also dominates the identity
(almost everywhere on real lines by differentiation).

\paragraph{Policies and paths.}
On a finite grid, an \emph{averaging policy} along axis $i$ is the
positive linear operator that chooses at each cell an interval along
axis $i$ through that cell and averages the input over it. It is
pointwise dominated by $M_i$ on nonnegative inputs. For a given
nonnegative input one can choose such a policy that attains $M_i$ at
each cell. Choosing these policies from right to left on the successive
maximal outputs gives a linear product attaining a prescribed word on
that input. A product $P_1\cdots P_k$ of policies is a sum over
\emph{paths} $(z_0,\ldots,z_k)$, where $z_{a}$ lies in the interval that
$P_a$ chooses at $z_{a-1}$; a path thus follows the averaging choices
from the output cell to the input cell, in the written order of the
word. When $S>0$ we write
\[
 \Phi(w)=\pair{\one_E}{M_w\one_F}/S ,
\]
suppressing the dependence on $E,F$ and the model.

\section{A planar crossing on the fixed dyadic trees}\label{planesection}
Work on $X\times Y$ with one dyadic tree on each axis. Let $E,F$ be
nonempty. Testing the root rectangle gives $0<\theta\le1$.
We prove that
\begin{equation}\label{planarbound}
 \Phi(12)\le C\theta^{1/6}(1+\log(1/\theta)).
\end{equation}
All rectangles tested in the argument are products of nodes of the two
trees that define the maxima. This point will allow us to use the
estimate on every transverse slice of a higher-dimensional grid.

The two maxima give horizontal segments rich in $E$ and vertical
segments rich in $F$. In each column $u$, choose disjoint dyadic
intervals $V\subset Y$ with
\[
 |F_u\cap V|\ge t_F |V|,\qquad 0<t_F\le1.
\]
Here $F_u=\{y:(u,y)\in F\}$. Call each chosen $\{u\}\times V$ a
marked vertical segment and write $B$ for their union. For a vertical
dyadic interval $K$, put
\[
 U_K=\{u:\text{a marked interval in column }u\text{ contains }K\},
 \qquad D_K(\rho)=\{u:M_x^{\D}\one_{U_K}(u)\ge\rho\}
\]
with $0<\rho\le1$. Here $M_x^{\D}$ acts on the line $X$.
The condition $u\in D_K(\rho)$ allows column $u$ to use $K$ as a comparison
interval whenever some horizontal dyadic interval through $u$ has at least
fraction $\rho$ of its columns in $U_K$.

Symmetrically, in each row $y$ choose disjoint horizontal dyadic
intervals $W\subset X$ with $|E^y\cap W|\ge t_E|W|$, where
$E^y=\{x:(x,y)\in E\}$ and $0<t_E\le1$; call each $W\times\{y\}$ a
marked horizontal segment and write $A$ for their union. For a horizontal
dyadic interval $H$ put
\[
 U'_H=\{y:\text{a marked interval in row }y\text{ contains }H\},
 \qquad D'_H(\rho)=\{y:M_y^{\D}\one_{U'_H}(y)\ge\rho\}.
\]

To use a common vertical interval at a crossing, we discard segments
whose $F$-density becomes too small on a comparison interval. The
following lemma bounds the size of the union of the discarded segments.

\begin{lemma}\label{remove}
For $\tau>0$, call a marked vertical segment $\{u\}\times V$ bad if there
is a dyadic interval $K\supseteq V$ with $u\in D_K(\rho)$ and
$|F_u\cap K|<\tau|K|$.
The union $\mathcal N_F$ of the bad segments satisfies
\[
 |\mathcal N_F|\le \frac{\tau}{t_F\rho}|B|.
\]
In particular every surviving segment retains density at least $\tau$
on \emph{every} comparison interval $K$ containing it. The same holds
for marked horizontal segments, with $E,t_E,A,U'_H,D'_H(\rho)$ in place of
$F,t_F,B,U_K,D_K(\rho)$.
\end{lemma}
\begin{proof}
Write $B^y=\{v:(v,y)\in B\}$ for the horizontal section of $B$. For
$y\in K$, $U_K\subseteq B^y$. Consequently
\[
 \mathcal O=\bigcup_K(D_K(\rho)\times K)
 \subseteq\{(u,y): M_x^{\D}\one_{B^y}(u)\ge\rho\},
 \qquad |\mathcal O|\le |B|/\rho .
\]
The last inequality follows by taking the maximal dyadic intervals
with density at least $\rho$ in each row. This estimates the union, even
when a single long segment appears in $U_K$ on many successive scales.

In column $u$, take the inclusion-maximal intervals $K$ that satisfy
both $u\in D_K(\rho)$ and $|F_u\cap K|<\tau|K|$. They are disjoint, lie in
the vertical section $\mathcal O_u=\{y:(u,y)\in\mathcal O\}$, and
contain the intervals of all bad segments in that column. In each such
$K$, disjointness of the original marked segments gives
\[
 \sum_{V\subseteq K\ {\rm bad}} |V|
 \le t_F^{-1}|F_u\cap K|\le (\tau/t_F)|K|.
\]
Sum over the disjoint maximal intervals in each column, then over the
columns, and use $|\mathcal O|\le |B|/\rho$. The horizontal statement is
the same argument with rows and columns interchanged.
\end{proof}

\subsection{A common rectangle at a crossing}
To bound the intersection of the horizontal and vertical maximal level
sets, we show that a large overlap of the two segment families forces
both sets to have substantial density on a single rectangle.

\begin{lemma}\label{cross}
Let $A$ be a union of marked horizontal segments of $E$-density at least
$t_E$ and $B$ a union of marked vertical segments of $F$-density at least
$t_F$, as above. Then
\[
 |A\cap B|\le 8\left(\frac{\theta}{t_E t_F}\right)^{1/6}
                     \sqrt{|A||B|}.
\]
\end{lemma}
\begin{proof}
If $A\cap B$ is empty the estimate is immediate. Otherwise set
$s=|A\cap B|$ and choose
\[
 \rho_x=\frac{s}{8|A|},\qquad \rho_y=\frac{s}{8|B|},\qquad
 \tau_F=t_F\rho_x\rho_y,\qquad \tau_E=t_E\rho_x\rho_y.
\]
Apply Lemma~\ref{remove} to the vertical segments with $\rho=\rho_x$ and
$\tau=\tau_F$. The bad ones cover at most $\rho_y|B|=s/8$ cells. The
horizontal form of Lemma~\ref{remove}, with $\rho=\rho_y$ and
$\tau=\tau_E$, removes horizontal segments covering at most
$\rho_x|A|=s/8$ cells. Let $G$ be the part of $A\cap B$ that survives
both deletions, so $|G|\ge3s/4$. The comparison intervals in both
applications of Lemma~\ref{remove} are determined by the original
families.

At each crossing $z=(u,y)\in G$, let $H_z$ and $V_z$ be the intervals
of its horizontal and vertical segments. To use $H_z\times V_z$ as a
test rectangle, we need enough columns with vertical intervals containing
$V_z$ to make $V_z$ a comparison interval, and enough columns with
surviving intervals inside $V_z$ to supply its $F$ density. We need the
same two counts with rows and columns exchanged for the $E$ density.

In each horizontal segment with interval $H$, sort the lengths
$|V_z|$ of the crossings in $G$ on that segment. Keep a crossing if
at least $\rho_x|H|$ entries are at most its length and at least
$\rho_x|H|$ are at least its length, counting equality in both tests.
At sorted position $r$, at least $r$ entries are no greater than
the entry there. Failure of the first test therefore requires
$r<\rho_x|H|$. Counting from the other end gives the same bound
for the second test. Thus at most $2\rho_x|A|=s/4$ crossings fail
over all horizontal segments. Apply the same rule to
the horizontal lengths on each vertical segment, with $\rho_y$ in
place of $\rho_x$. That cost is at most $2\rho_y|B|=s/4$. Both rules
count against the crossings in $G$, so at least $s/4$ pass both tests.

Choose such a crossing and put $H=H_z$, $K=V_z$.
Figure~\ref{fig:crossing} shows how the marks containing and contained
in the chosen intervals are used in each direction.

\begin{figure}[!ht]
\centering
\begin{tikzpicture}[
  x=0.9cm, y=0.9cm, font=\small,
  common rectangle/.style={draw=black!60, fill=black!7, line width=0.35pt},
  crossing segment/.style={black!60, line width=0.45pt},
  containing mark/.style={blue!65!black, densely dashed, line width=0.9pt},
  surviving mark/.style={orange!75!black, line width=1pt}
]
\node[font=\small\bfseries] at (2.65,4.8) {(a) Vertical marks for $F$};
\node[font=\small\bfseries] at (9.8,4.8) {(b) Horizontal marks for $E$};

\begin{scope}[shift={(1.25,0)}]
  \draw[common rectangle] (0,0) rectangle (2.8,2.2);
  \draw[crossing segment] (0,0.82) -- (2.8,0.82);
  \foreach \x in {0.45,0.82}
    \draw[containing mark] (\x,0) -- (\x,4.4);
  \foreach \x in {1.91,2.38}
    \draw[surviving mark] (\x,0) -- (\x,1.1);
  \fill (1.12,0.82) circle (1.5pt) node[above right] {$z$};
  \node[below] at (1.4,0) {$H$};
  \node[left] at (0,1.1) {$K$};
\end{scope}

\begin{scope}[shift={(7,1.1)}]
  \draw[common rectangle] (0,0) rectangle (2.8,2.2);
  \draw[crossing segment] (1.12,0) -- (1.12,2.2);
  \foreach \y in {0.32,1.78}
    \draw[containing mark] (0,\y) -- (5.6,\y);
  \foreach \y in {0.57,1.4}
    \draw[surviving mark] (0,\y) -- (1.4,\y);
  \fill (1.12,0.82) circle (1.5pt) node[above right] {$z$};
  \node[below] at (1.4,0) {$H$};
  \node[left] at (0,1.1) {$K$};
\end{scope}

\draw[containing mark] (1.25,-0.8) -- (1.85,-0.8)
  node[right] {Containing mark};
\draw[surviving mark] (6.25,-0.8) -- (6.85,-0.8)
  node[right] {Contained surviving mark};
\end{tikzpicture}
\caption{The shaded rectangle is $H\times K$ in both panels. In (a),
marks containing $K$ certify $H\subseteq D_K(\rho_x)$, and surviving marks
inside $K$ supply the $F$-density bound. In (b), marks containing $H$
certify $K\subseteq D'_H(\rho_y)$, and surviving marks inside $H$ supply
the $E$-density bound.}
\label{fig:crossing}
\end{figure}
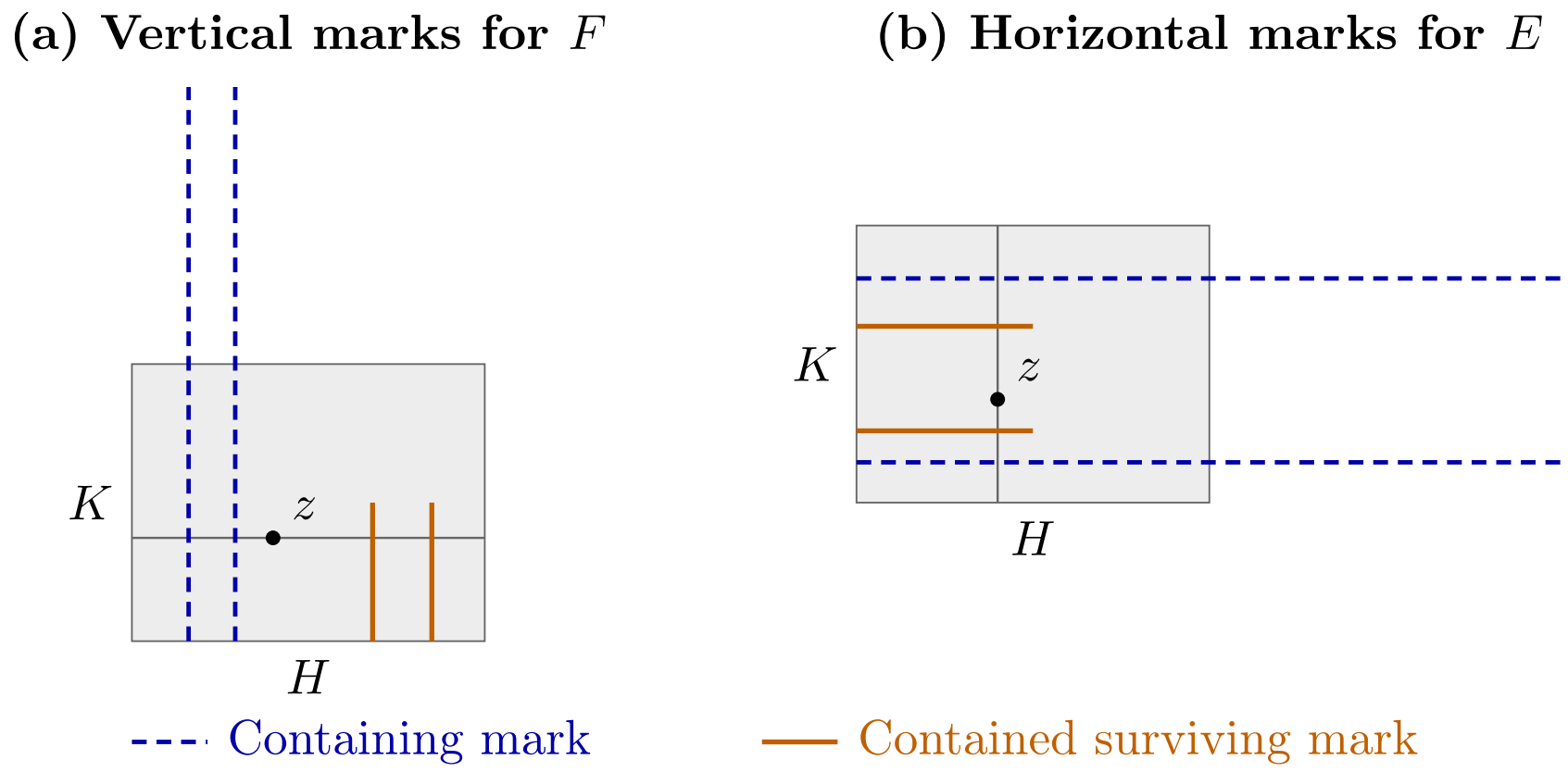

At least $\rho_x|H|$
columns along the row segment have a crossing in $G$ whose vertical
interval has length at least $|K|$. That interval contains $K$ by dyadic
nesting through $y$. Thus $H$ has $U_K$ density at least $\rho_x$, and
\emph{all} $u'\in H$ belong to $D_K(\rho_x)$. At least $\rho_x|H|$
columns have a crossing in $G$ with vertical interval contained in $K$.
Each such segment survived removal, so $|F_{u'}\cap K|/|K|\ge\tau_F$.
Averaging in $u'$ yields $\av{\one_F}{H\times K}\ge\rho_x\tau_F$.

In the other direction, at least $\rho_y|K|$ rows have a horizontal
interval containing $H$, so $K\subseteq D'_H(\rho_y)$. At least
$\rho_y|K|$ rows have a surviving interval contained in $H$. This yields
$\av{\one_E}{H\times K}\ge\rho_y\tau_E$.
Testing the same rectangle with both sets gives
\[
 \theta\ge t_Et_F(\rho_x\rho_y)^3
       =t_Et_F\left(\frac{s^2}{64|A||B|}\right)^3 .
\]
Solving for $s$ gives the bound.
\end{proof}

\subsection{From crossings to the maximal operators}
Use maximal disjoint dyadic intervals at each level, so
$A_q=\{M_x^{\D}\one_E>q\}$ and $B_t=\{M_y^{\D}\one_F>t\}$ are unions
of marked segments with $t_E=q$ and $t_F=t$. For $0<q,t<1$,
Lemma~\ref{cross} and the one-dimensional weak $(1,1)$ bounds
$|A_q|\le|E|/q$, $|B_t|\le|F|/t$ give
\begin{align}
 |A_q\cap B_t|&\le 8\theta^{1/6}(qt)^{-2/3}\sqrt{|E||F|},\label{levels}\\
 \sum_{B_t} M_x^{\D}\one_E&\le
       24\theta^{1/6}t^{-2/3}\sqrt{|E||F|}.\label{integral}
\end{align}
The second inequality integrates the first in $q$. To bound the
composition, we still need to estimate what the outer horizontal
maximal operator does to each $\one_{B_t}$.

For any set $B$ and $0<v<1$, in each row take the maximal intervals
with $B$ density greater than $v$. For each such interval $I$, all
its cells have $M_x^{\D}\one_E$ at least the row density of $E$ on $I$.
Summing on their disjoint $B$ portions gives
\[
 |E\cap\{M_x^{\D}\one_B>v\}|\le v^{-1}\sum_B M_x^{\D}\one_E.
\]
Split the layer integral at $0<\epsilon<1$. For the part below
$\epsilon$, the weak $(1,1)$ bound gives
$\|\min(\epsilon,M_x^{\D}\one_B)\|_2^2\le2\epsilon|B|$
by integrating $2v|\{M_x^{\D}\one_B>v\}|$. For the part above
$\epsilon$, integrate the preceding estimate in $v$. Together with
Cauchy--Schwarz for the first part, this gives
\[
 \sum_E M_x^{\D}\one_B
 \le \sqrt{2\epsilon |E||B|}
          +\log(1/\epsilon)\sum_B M_x^{\D}\one_E .
\]
Apply this with $B=B_t$. By the layer-cake formula and sublinearity,
\[
 M_x^{\D}M_y^{\D}\one_F\le\int_0^1 M_x^{\D}\one_{B_t}\,dt .
\]
Since $|B_t|\le |F|/t$, integrals of $t^{-1/2}$ and $t^{-2/3}$ are
finite. It follows that
\[
 \frac{\sum_E M_x^{\D}M_y^{\D}\one_F}{\sqrt{|E||F|}}
 \le 2\sqrt{2\epsilon}+72\theta^{1/6}\log(1/\epsilon).
\]
For $0<\theta<1$, take $\epsilon=\theta^{1/3}$ to obtain
\eqref{planarbound}. For $\theta=1$ the $L^2$ bound suffices. It follows
that $\Phi(12)\le C_\gamma\theta^\gamma$ for every $0<\gamma<1/6$.
Interchanging the axes gives the other order.

\section{Adding a coordinate}\label{dimension}
To extend the estimate to a new coordinate, we apply the transverse word
at each height. The next argument bounds the slice density products by
testing rectangles in the full grid after deleting small exceptional
sets.
\begin{proposition}\label{add}
Let $X$ be a dyadic product grid, let $T$ be a product of coordinate
maximal operators visiting each axis of $X$ once, and suppose that for
some $\gamma>0$ and $C_T<\infty$, uniform in the grid size,
\[
 \pair{\one_E}{T\one_F}\le C_T\,\theta_X(E,F)^\gamma\sqrt{|E||F|}
 \qquad\text{for all }E,F\subseteq X,
\]
where $\theta_X$ is taken over dyadic rectangles of $X$. Let $Z$ be
another dyadic line and fix $1<p<2$. On $X\times Z$ the word $TM_Z$ has
a grid-independent estimate of the same form with exponent
\begin{equation}\label{rec}
 \gamma'=\bigl(6+(2p+1)/\gamma\bigr)^{-1}.
\end{equation}
\end{proposition}
\begin{proof}
Assume $E,F$ are nonempty. For $0<t<1$ put $B=\{M_Z\one_F>t\}$ and
$\delta=(\theta/t)^{\gamma'}$. We will bound the pairing with $T\one_B$
by $C\delta\sqrt{|E||B|}$ and then integrate in $t$. If $B$ is empty,
the level estimate is immediate. Otherwise mark the maximal intervals
with $F$ density greater than $t$ in each vertical fiber. They are
disjoint there, and $|B|\le |F|/t$. For a vertical node $K$ let
$U_K\subset X$ be the bases with a marked interval containing $K$.
For $0<\rho<1$, let $D_K=\{M_X\one_{U_K}\ge\rho\}$, where $M_X$
is the maximal operator over all transverse dyadic rectangles. At
height $s\in K$, $U_K$ is contained in the slice $B_s$. The union of
$D_K$ over such $K$ is therefore contained in
$\{M_X\one_{B_s}\ge\rho\}$. Successive coordinate maxima dominate
$M_X$ on nonnegative inputs, so $M_X$ is bounded on $L^p(X)$. Applying
Chebyshev's inequality and the $L^p$ bound on each slice gives
$|\bigcup_K D_K\times K|\le C_{d,p}\rho^{-p}|B|$.

Choose density thresholds $u_E,u_F>0$, to be fixed below. At each
base $x$, take the maximal nodes $K$ satisfying both $x\in D_K$
and $|E^x\cap K|<u_E|K|$, where superscript $x$ denotes a vertical
section. Delete the $E$ points in those nodes. For $F$, take the maximal
nodes satisfying $x\in D_K$ and $|F^x\cap K|<u_F|K|$, and delete from
$B$ the marks contained in them. In each such node the deleted marks
have total length at most $|F^x\cap K|/t$. In each fiber over $x$,
the maximal nodes for each deletion are disjoint, and their segments
$\{x\}\times K$ lie in $\bigcup_K(D_K\times K)$. The surviving sets
$E',B'$ therefore satisfy
\[
 |E\setminus E'|\le C_{d,p}u_E\rho^{-p}|B|,\qquad
 |B\setminus B'|\le C_{d,p}(u_F/t)\rho^{-p}|B|.
\]
In a transverse rectangle $H$ at height $s$, let $\pi,b$ be the
densities of the surviving slices $E'_s,B'_s$ on $H$. If $b\ge2\rho$,
sort the marked intervals at the bases of $B'_s\cap H$ by length and
choose a middle entry $K$, allowing either middle entry for an even
count. The intervals all contain $s$ and hence are nested.
At least half contain $K$, so the $U_K$ density on $H$ is at least
$b/2\ge\rho$ and $H\subset D_K$. At least half are contained in $K$.
At any base in $H$, if the vertical $E$ density on $K$ were below
$u_E$, then $K$ would lie inside a maximal node used for the $E$
deletion and any $E$ point at height $s$ would have been deleted. The
same reasoning for $F$ shows that a surviving mark inside $K$ requires
vertical $F$ density at least $u_F$ on $K$ at its base. Thus the $E$
and $F$ densities on $H\times K$ are at least $\pi u_E$ and $bu_F/2$,
respectively, so $\pi bu_Eu_F/2\le\theta$.
If $b<2\rho$, use $\pi b\le b$. In either case the slice density product
is at most $2\max(\rho,\theta/(u_Eu_F))$.

We choose the thresholds to make both the deletion error and the
surviving pairing at most a constant times
$\delta\sqrt{|E||B|}$. The $L^2$ bound of $T$ gives this estimate
directly if $\delta\ge1$. Preservation of constants gives
$T\one_B\le1$, hence $\pair{\one_E}{T\one_B}/\sqrt{|E||B|}
\le (|E|/|B|)^{1/2}$. This is at most $\delta$ if $|B|\ge\delta^{-2}|E|$.
For $0<\delta<1$ and $0<|B|<\delta^{-2}|E|$, choose
$\rho=\delta^{1/\gamma}$ so that $\rho^\gamma=\delta$, and set
$u_E=\delta^4\rho^p$ and $u_F=t\delta^2\rho^p$. The recurrence
\eqref{rec} ensures
$u_Eu_F\rho=t\delta^6\rho^{2p+1}=t\delta^{6+(2p+1)/\gamma}=\theta$.
The deletion costs are at most $C_{d,p}\delta^4|B|\le C_{d,p}\delta^2|E|$
and $C_{d,p}\delta^2|B|$, respectively. Sublinearity bounds the errors by
the pairings for $(E\setminus E',B)$ and $(E',B\setminus B')$. Their
sum is at most $C\delta\sqrt{|E||B|}$ by the $L^2$ bound.
The slice density product is at most $2\rho$, so the transverse estimate
and Cauchy--Schwarz in the height $s$ give
\[
 \pair{\one_{E'}}{T\one_{B'}}
 \le C_T(2\rho)^\gamma\sum_s\sqrt{|E'_s||B'_s|}
 \le 2^\gamma C_T\delta\sqrt{|E||B|}.
\]
To pass from the level estimate to the full pairing, write
$M_Z\one_F$ as the integral of its level indicators. This integral
is a finite sum on the grid, so sublinearity gives
$TM_Z\one_F\le\int_0^1 T\one_{\{M_Z\one_F>t\}}\,dt$. Since
$|B|^{1/2}\le(|F|/t)^{1/2}$, the integrated pairing is at most
$C\theta^{\gamma'}S\int_0^1t^{-1/2-\gamma'}\,dt$,
which is finite because $\gamma'<1/2$.
\end{proof}
The map $\gamma\mapsto(6+3/\gamma)^{-1}$ is increasing and continuous,
and the exponent \eqref{rec} increases to $(6+3/\gamma)^{-1}$ as
$p\downarrow1$. Choose the planar exponent $\gamma_2$ in
\eqref{planarbound} arbitrarily close to $g_2=1/6$ from below and take
$p>1$ sufficiently close to $1$ at each application of
Proposition~\ref{add}. By continuity this yields estimates for a word
visiting each axis once with exponents arbitrarily close to $g_d$
from below, where $g_d$ is as in Theorem~\ref{word}. Since $\theta\le1$,
these estimates give every exponent $0<\gamma_d<g_d$.
Permuting the axes gives all orders. In dimension three this is every
$0<\gamma_3<1/24$.

\section{Reducing repeated coordinates}\label{compress}
To obtain a density exponent independent of the word length, we
remove all but one visit to a chosen coordinate in a single reduction.

Fix a word $w=i_1\cdots i_k$ on a dyadic product grid, with nonempty
test sets $E,F$; if either set is empty, all pairings vanish. Choose a
coordinate $j$ that occurs in $w$. For each position $t$ with $i_t=j$,
split $w$ into the prefix before $t$, the step $j$, and the suffix after
$t$. Write $u_t,v_t$ for the prefix and suffix with every occurrence of
$j$ deleted, and $v^\rev$ for the reversal of a word $v$.

\begin{lemma}\label{reduce}
For $0<\lambda<1$ one has
\[
 \Phi(w)\le k C_\lambda^k
       \max_{t:\,i_t=j}\Phi(u_t j v_t)^{\lambda^2}.
\]
\end{lemma}
\begin{proof}
The proof has four steps.
\begin{enumerate}
\item Stopping policies $P_a$ with $M_w\one_F\le2^kP_1\cdots P_k\one_F$,
 and a factorization \eqref{kern} of each $P_a$ into a forward kernel
 $L_a$ bounded on $L^2$ and a reverse kernel $R_a$ dominated by a
 maximal operator.
\item Tagging each path by a longest $j$-interval $J$ and applying
 H\"older to the suffix after it:
 $m_t\le C_\lambda^kS^{1-\lambda}(m^{\rm rev}_t)^\lambda$.
\item In the reverse sum, keeping only the pivot $j$-step:
 $m^{\rm rev}_t\le\pair{\one_F}{M_{v_t^\rev}g}$ with
 $g=A^*\one_E$ and $A\le M_{u_t}M_j$.
\item Repeating steps (1)--(2) on $M_{v_t^\rev}$ to restore the order:
 $m^{\rm rev}_t\le C_\lambda^kS^{1-\lambda}
 \pair{\one_E}{M_{u_t}M_jM_{v_t}\one_F}^\lambda$.
\end{enumerate}

\paragraph{Stopping intervals and reverse kernels.}
At step $a$ use the nonnegative input
$g_a=M_{i_{a+1}}\cdots M_{i_k}\one_F$, with $g_k=\one_F$. On each fiber
along axis $i_a$, retain the root. For each retained interval $I$, let
$\mathcal J(I)$ consist of the inclusion-maximal members among the
strict dyadic subintervals $J\subsetneq I$ satisfying
\[
 \av{g_a}J>2\av{g_a}I.
\]
Retain the members of $\mathcal J(I)$. Their total length is at most
$|I|/2$.
Choose the smallest retained interval $I_z$ through $z$. The policy
$P_a$ averages on $I_z$. Each dyadic interval through $z$ has average at
most twice that of its smallest retained ancestor by the stopping rule,
and retained averages increase down the chain through $z$. Thus
$P_ag_a\ge\frac12M_{i_a}g_a$. Positivity bounds the given word by
$2^kP_1\cdots P_k$ on $\one_F$.

For a retained interval $I$, let $n_I(z')$ count its strict retained
descendants through $z'\in I$. The descendants at generation $\ell$
below $I$ cover $\{z'\in I:n_I(z')\ge\ell\}$ and have total length
at most $2^{-\ell}|I|$. Suppress the fixed transverse coordinates of
the fiber, put
$\epsilon=(1-\lambda)/(8\lambda)$, $c=1-2^{-\epsilon}$, and define
\[
 L_a(z,z')=P_a(z,z')2^{n_{I_z}(z')/8}\quad(z'\in I_z),\qquad
 R_a(z',z)=c\sum_{I\ {\rm retained}:\,z',z\in I}
                  \frac{2^{-\epsilon n_I(z')}}{|I|}.
\]
Take $L_a=0$ off the forward support. The count $n_I(z')$ takes the
successive values $0,1,\ldots$ as $I$ runs from the smallest retained
interval containing $z'$ to the root. Thus $R_a$ is a mixture of
averages through its departure point $z'$ with total mass at most
$c\sum_{n\ge0}2^{-\epsilon n}=1$. On nonnegative inputs it is
pointwise dominated by $M_{i_a}$. The tail bound gives
$\av{2^{3 n_I/8}}I\le\sum_{\ell\ge0}2^{-5\ell/8}$.
H\"older with exponents $3,3/2$ therefore gives
$L_ah\le C(M_{i_a}h^{3/2})^{2/3}$ for $h\ge0$, hence
$\|L_a\|_{2\to2}\le C$. Keeping only the term $I=I_z$ in $R_a$ cancels
the depth powers, since $\epsilon\lambda=(1-\lambda)/8$, and gives
\begin{equation}\label{kern}
 P_a(z,z')\le c^{-\lambda}L_a(z,z')^{1-\lambda}R_a(z',z)^\lambda .
\end{equation}
All three kernels preserve the coordinates other than $i_a$.

\paragraph{Choosing the longest interval and reversing the suffix.}
Write the policy pairing as the sum over paths $(z_0,\ldots,z_k)$ of
\[
 \one_E(z_0)\prod_{a=1}^k P_a(z_{a-1},z_a)\,\one_F(z_k).
\]
Tag each path of positive weight by the earliest time $t$ of a longest
chosen $j$-interval $J$. Between two successive $j$-steps the
$j$-coordinate stays fixed, so their intervals intersect. Inductively
the union up to any such step is one of its intervals. The next
interval meets that union, and intersecting nodes in the common tree
are nested. Thus $J$ contains every visited $j$-interval and both
endpoint $j$-coordinates.
Let $m_t$ be the pairing over paths with tag $t$, and on those paths
put $b=\one_E(z_0)\one_F(z_k)\prod_{a\le t}P_a(z_{a-1},z_a)$; this
includes the pivot factor $P_t(z_{t-1},z_t)=|J|^{-1}$. Applying
\eqref{kern} to the suffix and H\"older's inequality with exponents
$1/(1-\lambda)$ and $1/\lambda$ gives
\[
 m_t\le c^{-\lambda k}
   \Bigl(\sum_{\text{tag }t}b\prod_{a>t}L_a(z_{a-1},z_a)\Bigr)^{1-\lambda}
   \Bigl(\sum_{\text{tag }t}b\prod_{a>t}R_a(z_a,z_{a-1})\Bigr)^{\lambda}.
\]
Dropping the path restrictions in the first sum bounds it by
$\pair{\one_E}{P_1\cdots P_tL_{t+1}\cdots L_k\one_F}\le C^kS$.
Writing $m^{\rm rev}_t$ for the second sum gives
\[
 m_t\le C_\lambda^k S^{1-\lambda}(m^{\rm rev}_t)^\lambda .
\]

\paragraph{Removing the other steps in coordinate $j$.}
Write a cell as $(x,y)$, with transverse vector $x$ and $j$-coordinate
$y$. Fix a prefix ending at the cell $z_{t-1}=(x,y)$ just before the
selected step, and let $y_0$ be the starting $j$-coordinate. This last
prefix cell determines the interval $J$. The selected step arrives at
$z_t=(x,y')$ and the complete path ends at $z_k=(\xi,s)$.
The reverse product $Q=R_k\cdots R_{t+1}$ connects $(\xi,s)$ to
$(x,y')$ and is the identity if the suffix is empty. Summing over
its arrival coordinate $y'$ defines a kernel on each endpoint slice:
\[
 \Delta_s(\xi,x)=\sum_{y'} Q((\xi,s),(x,y')),\qquad
 \mathcal V((\xi,s),(x,s''))=\one_{\{s=s''\}}\Delta_s(\xi,x).
\]
On a tagged continuation, $s,y_0,y'$ are all in $J$ and the selected
step contributes $|J|^{-1}$. Keep the conditions $s,y_0\in J$ and that
single factor, then drop the restriction on $y'$ and the remaining
restrictions on the suffix; positivity lets us sum over all $y'$. This
bounds the sum of tagged continuations, without the prefix weight, by
\begin{equation}\label{pivot}
 \frac{\one_{\{y_0\in J\}}}{|J|}\sum_{s\in J}\sum_\xi
     \Delta_s(\xi,x)\one_F(\xi,s)
 =\one_{\{y_0\in J\}}\av{(\mathcal V^*\one_F)(x,\cdot)}J
 \le (M_j\mathcal V^*\one_F)(x,y_0).
\end{equation}
Here $(\mathcal V^*\one_F)(x,s)=\sum_\xi\Delta_s(\xi,x)\one_F(\xi,s)$ by
transposition in counting measure. The final inequality uses the
fact that $J$ contains $y_0$ whenever the indicator equals one.

To bound $\mathcal V$ on nonnegative $h$, extend the slice at $s$
constantly in the running $j$-coordinate $y'$. By the definition of
$\Delta_s$,
\[
 h^{[s]}(x,y'):=h(x,s),\qquad
 (Qh^{[s]})(\xi,s)
 =\sum_{x,y'}Q((\xi,s),(x,y'))\,h(x,s)=(\mathcal Vh)(\xi,s).
\]
Evaluate $Q$ from right to left on $h^{[s]}$. Each $j$-step has row
mass at most one and preserves any nonnegative upper bound independent
of the running $j$-coordinate. Each transverse step applies at most
its own coordinate maximal operator to that bound. Since
$Q=R_k\cdots R_{t+1}$, the resulting bound at $(\xi,s)$ is
$(M_{v_t^\rev}h)(\xi,s)$.

Now fix the starting cell $(x_0,y_0)$ of the forward prefix. With
$\Psi=M_j\mathcal V^*\one_F$, the right side of \eqref{pivot} bounds the
continuations at the last prefix cell $(x,y)$ by
\[
 \Psi^{[y_0]}(x,y)=\Psi(x,y_0)
 \qquad\text{for every }y.
\]
Sum this bound over prefixes from $(x_0,y_0)$. Reading the prefix
product from right to left, each $j$-policy preserves the upper bound
independent of the running $j$-coordinate, and each transverse
policy applies at most its coordinate maximal operator to it. This
leaves the transverse maxima in the written order $u_t$, so
\[
 (P_1\cdots P_{t-1}\Psi^{[y_0]})(x_0,y_0)
 \le (M_{u_t}\Psi)(x_0,y_0).
\]
Summing over starting cells in $E$ yields
$m^{\rm rev}_t\le\pair{\one_E}{M_{u_t}M_j \mathcal V^*\one_F}$. Choose a
linear averaging product $A$ attaining $M_{u_t}M_j$ on
$\mathcal V^*\one_F$. It satisfies $A\le M_{u_t}M_j$ on nonnegative
inputs. With $g=A^*\one_E$ this yields
\begin{equation}\label{second}
 m^{\rm rev}_t\le\pair{\one_F}{\mathcal Vg}\le\pair{\one_F}{M_{v_t^\rev}g},
 \qquad \|g\|_2\le C^k\sqrt{|E|}.
\end{equation}

\paragraph{Restoring the suffix order.}
Apply the stopping construction to $M_{v_t^\rev}$ on input $g$, then
use \eqref{kern} on all its letters and apply H\"older. The pairing
with the forward kernels is at most
$C^k\|\one_F\|_2\|g\|_2\le C^kS$, after enlarging $C$.
If $R'_1,\ldots,R'_{\ell}$ are the new reverse kernels in the written
order of $v_t^\rev$, the reverse pairing is
\[
 \pair{\one_F}{(R'_1)^*\cdots(R'_{\ell})^*g}
 =\pair{g}{R'_{\ell}\cdots R'_1\one_F}
 \le\pair{g}{M_{v_t}\one_F}
 \le\pair{\one_E}{M_{u_t}M_j M_{v_t}\one_F}.
\]
The transposition restores the original order of the transverse suffix.
If $v_t$ is empty, both pairings reduce to $\pair{\one_F}{g}$,
whose powers $1-\lambda$ and $\lambda$ satisfy the same bounds.
Thus
\[
 m^{\rm rev}_t\le C_\lambda^k S^{1-\lambda}
            \pair{\one_E}{M_{u_t}M_jM_{v_t}\one_F}^{\lambda}.
\]
Substituting in the bound for $m_t$ gives
$m_t/S\le C_\lambda^k\Phi(u_t j v_t)^{\lambda^2}$.
Sum over at most $k$ tags and include the initial factor $2^k$.
\end{proof}

\subsection{Iteration of the reduction}
At each reduction, choose a position attaining the maximum in
Lemma~\ref{reduce}. The replacement removes all but one occurrence of
the chosen coordinate and preserves the multiplicities and relative
order of the others. Thus after $r\le d$ reductions every remaining
coordinate occurs once. If some axes are absent, append their maximal
operators on the right. These dominate the identity on nonnegative
inputs, so the padded word gives an upper bound and has the estimate of
Sections~\ref{planesection}--\ref{dimension} for the full $d$-dimensional
density $\theta$.

At each reduction the word length stays at most $k$, so the
prefactor is bounded by $(2C_\lambda)^k$. There are at most $d$ such
factors. Raising them to exponents at most one can only decrease
them when the constants are chosen at least one. The resulting density
exponent is $\lambda^{2r}\gamma_d$. Given $0<\eta_d<g_d$, choose
$\eta_d<\gamma_d<g_d$ as in Section~\ref{dimension} and $\lambda<1$ with
$\lambda^{2d}\gamma_d\ge\eta_d$. Since $\theta\le1$, the exponent
$\eta_d$ is valid. This proves the dyadic case of Theorem~\ref{word}.

\section{Ordinary intervals and Lebesgue measure}\label{ordinary}
\subsection{Ordinary finite grids}
We compare the ordinary averaging paths with paths on randomly shifted
dyadic trees. Choose maximizing ordinary policies for the word on input
$\one_F$. Pad each axis of $n$ cells to a cycle of dyadic length
$N\ge4n$, keeping the original cells in one arc and adding no points
to $E$ or $F$. On each axis shift the binary tree by a uniformly chosen
cell offset from $0$ to $N-1$. Read the cycle as a line starting at
the shifted image of the first leaf in the unshifted tree. This gives
a dyadic interval tree in the resulting order. All fibers and all visits
to that axis use this one tree, as required by the longest-interval
reduction.

For a chosen interval $I$ of length $\ell$, let $\widehat I$ be its
smallest containing node in the shifted tree, allowing the root. If $b$
is the least dyadic length at least $\ell$, the probability that $I$
crosses a boundary at scale $2^m b<N$ is $(\ell-1)/(2^m b)\le2^{-m}$.
The logarithmic dilation $\log_2(|\widehat I|/\ell)$ consists of a
rounding term $\log_2(b/\ell)<1$ plus one unit for each such scale at
which $I$ crosses a boundary. Summing these tail probabilities gives
$\mathbb E\log_2(|\widehat I|/\ell)\le3$. Replacing each interval by
$\widehat I$ preserves every original path. For a fixed path with
intervals $I_1,\ldots,I_k$, Jensen's inequality gives its expected
weight ratio
\[
 \mathbb E\,2^{-\sum_{a=1}^k \log_2(|\widehat I_a|/|I_a|)}
 \ge 2^{-\sum_{a=1}^k\mathbb E\log_2(|\widehat I_a|/|I_a|)}\ge 8^{-k}.
\]
Jensen's inequality uses only the individual expected log dilations,
so repeated visits can use the same shift. Extend the policies to
padded cells by singleton averages. Extra paths have nonnegative weight.
Each enlarged policy product is pointwise dominated on nonnegative
inputs by the word on the corresponding shifted dyadic grid.

A proper node has length at most $N/2$. To meet the original arc in
two pieces it would have to contain the entire complementary arc,
of length at least $3N/4$. Thus each node meets the original arc in
a consecutive interval or the empty set, and the root contains the
whole arc. Intersecting any shifted dyadic rectangle with the
original grid preserves the counts of $E,F$ and can only increase
their densities when the intersection is nonempty. Thus the dyadic
density test is at most the ordinary one. Summing the expected weights
of the original paths bounds $8^{-k}\pair{\one_E}{M_w\one_F}$, with
ordinary maxima, by the expected pairing with shifted dyadic maxima.
Apply the dyadic estimate inside this expectation and use the density
comparison to obtain the ordinary-grid bound.

\subsection{The continuum word estimate}
For finite-measure $E,F\subset\mathbb R^d$, partition
$Q_n=[-2^n,2^n)^d$ into cells of side $2^{-n}$. Let $E_n$ be the union
of the cells with $E$ density above $1/2$, and define $F_n$ in the
same way. On $Q_n$, the pointwise error of each thresholded
indicator is at most twice the error of the corresponding dyadic
cell average. The cell averages converge in $L^2$, and the $L^2$
tails of $\one_E$ and $\one_F$ outside $Q_n$ tend to zero. Thus
$\one_{E_n}\to\one_E$ and $\one_{F_n}\to\one_F$ in $L^2$. On every
aligned rectangle the density of $E_n$ is at most twice that of $E$,
and likewise for $F_n$. The ordinary grid density test for the selected
cells is therefore at most $4\theta(E,F)$.

Let $\mathcal I_n$ consist of the nonempty unions of consecutive half-open
mesh intervals in $[-2^n,2^n)$. For $h\in L^2(\mathbb R^d)$ set, on
almost every fiber in $Q_n$,
\[
 M_{i,n}h(x)=
 \max_{\substack{I\in\mathcal I_n\\x_i\in I}}
 \frac{1}{|I|}\int_I |h(x_1,\ldots,x_{i-1},t,x_{i+1},\ldots,x_d)|\,dt .
\]
Set $M_{i,n}h=0$ outside $Q_n$. Extend a grid array $a$ to $\widetilde a$
by its value $a_m$ on cell $m$ and by zero off $Q_n$. On this cell the
fiber average of $|\widetilde a|$ over mesh intervals $r,\ldots,s$ is
\[
 \frac{2^{-n}\sum_{\ell=r}^s |a_{m_1,\ldots,m_{i-1},\ell,m_{i+1},\ldots,m_d}|}
      {(s-r+1)2^{-n}} .
\]
The intervals through $x_i$ correspond exactly to the consecutive index
sets containing $m_i$. Hence $M_{i,n}\widetilde a=\widetilde{M_i^{\rm grid}a}$
almost everywhere. Iterate this identity along the word. In passing
from counting measure to Lebesgue measure, multiply the pairing and the
geometric mean $\sqrt{\#E_n\,\#F_n}$ of the two cell counts by the cell
volume. The geometric mean becomes $\sqrt{|E_n||F_n|}$, giving the
aligned word estimate for the selected indicators.

Approximating interval endpoints on each locally integrable fiber gives
$M_{i,n}h\uparrow M_i h$ almost everywhere. Domination by $M_i h\in L^2$
gives $L^2$ convergence. If $h_n\to h$ in $L^2$, then
$|M_{i,n}h_n-M_{i,n}h|\le M_i(h_n-h)$ gives $M_{i,n}h_n\to M_i h$ in
$L^2$ as well. Iterating gives convergence of
$M_{i_1,n}\cdots M_{i_k,n}\one_{F_n}$ to $M_w\one_F$ in $L^2$.
Pair with $\one_{E_n}$ and pass to the limit. Absorbing $4^{\eta_d}$
into the base gives Theorem~\ref{word} on $\mathbb R^d$.

\section{Two-weight bounds on a finite grid}\label{weights}
On an ordinary finite grid, we use Theorem~\ref{word} to bound the
$L^2$ norm of a policy product restricted to a pair of cell sets.
Applying this bound to the level bands of two weights will give the
weighted norm estimate.

A policy may now choose both its coordinate and its interval at each
cell. By sublinearity, a product $P$ of $k$ such policies on nonnegative
inputs is bounded pointwise by the sum of all $d^k$ words of length $k$.
Thus the density estimate and unweighted $L^p$ bounds hold for $P$ with
at most exponential growth in $k$. Fix $0<\eta_d<g_d$ and
$0<\beta<\sigma<\eta_d$ depending only on $d$.

\begin{lemma}\label{localblock}
For a product $P$ of $k\ge1$ policies and cell sets $U,V$, use $\one_V$
to restrict the input of $P$ and $\one_U$ to restrict its output. Then
\begin{equation}\label{block}
 \|\one_U P\one_V\|_{2\to2}\le C_d^k\theta(U,V)^\sigma .
\end{equation}
\end{lemma}
\begin{proof}
Take $U'\subset U$, $V'\subset V$ with $|U'|=m$, $|V'|=n$.
The $L^4$ and $L^{4/3}$ bounds give
\[
 \pair{\one_{U'}}{P\one_{V'}}
 \le C_d^k\min(m^{1/4}n^{3/4},m^{3/4}n^{1/4}).
\]
Raise the density bound to the power $\sigma/\eta_d$ and the
unweighted bound to $1-\sigma/\eta_d$, then multiply. The density
test for the subsets is at most $\theta(U,V)$. Both interpolated
monomials increase in $m,n$. For $m\le2^i,n\le2^j$, evaluating them
at these upper bounds and taking the smaller gives
\[
 C_d^k\theta(U,V)^\sigma\,2^{(i+j)/2}
      2^{-(1-\sigma/\eta_d)|i-j|/4}
\]
for the indicator pairing.
For nonnegative tests $h,g$ supported in $U,V$, respectively, rank
the values of each test in decreasing order, starting at rank one.
Let $h_i,g_j$ be their maxima on the respective rank blocks
$[2^i,2^{i+1})$ and $[2^j,2^{j+1})$, with value zero for empty blocks.
For $i\ge1$ with a nonempty block, $2^ih_i^2$ is at most twice the
sum of squares on the full preceding block. Thus
$\sum_i2^ih_i^2\le3\|h\|_2^2$, and likewise for $g$.

Multiply the indicator bound on each pair of blocks by $h_i g_j$
and sum. After taking out $C_d^k\theta(U,V)^\sigma$, this is the
$\ell^2$ pairing of $(2^{i/2}h_i)_i$ with the convolution of
$(2^{j/2}g_j)_j$ and the summable kernel
$\kappa_m=2^{-(1-\sigma/\eta_d)|m|/4}$. Young's inequality gives
\eqref{block}. Positivity handles signed tests.
\end{proof}

The preceding lemma controls the operator on one pair of level
bands. Summing these bounds gives the two-weight estimate needed
to control the iteration in the next section.

\begin{proposition}\label{weighted}
Suppose $a,b$ are nonnegative real-valued grid functions with
$\av aR\av bR\le1$ on all ordinary grid rectangles. Let $D_a,D_b$ be
multiplication by $a^\beta,b^\beta$, respectively. For every product
$P$ of $k\ge0$ policies, $\|D_a P D_b\|_{2\to2}\le C_d^{k+1}$.
\end{proposition}
\begin{proof}
For $k\ge1$, let $U_i=\{2^i\le a<2^{i+1}\}$ and
$V_j=\{2^j\le b<2^{j+1}\}$ for $i,j\in\mathbb Z$. On every rectangle,
$\av{\one_{U_i}}R\le2^{-i}\av aR$ and
$\av{\one_{V_j}}R\le2^{-j}\av bR$. Thus
$\theta(U_i,V_j)\le\min(1,2^{-i-j})$. By Lemma~\ref{localblock},
the weighted block norm from $V_j$ to $U_i$ is at most
$C_d^{k+1}q_{i+j}$, where
\[
 q_m=2^{\beta m}\min(1,2^{-\sigma m}).
\]
Since $0<\beta<\sigma$, this kernel is summable over
$m\in\mathbb Z$. Sum the block bounds against the bandwise $L^2$
norms of the two tests. Reflecting $j$ to $-j$ makes this sum a
convolution estimate, so Young's inequality and disjointness of the
bands give $\|D_a P D_b\|_{2\to2}\le C_d^{k+1}$. Zero weights
contribute nothing. For the empty word use the singleton bound
$ab\le1$.
\end{proof}

\section{A finite barrier and a bounded correction}\label{barrier}
The weighted bounds give a grid function between a lower and an upper
logarithmic bound, with oscillation controlled independently of the
grid. For the exponential weights of a real $f$, subtract this function
from the cell averages of $\alpha f$ and divide by $\alpha$ to obtain
the bounded correction for the grid approximation of $f$.

\begin{lemma}\label{finitebarrier}
For the weights in Proposition~\ref{weighted} there is a real grid
function $G$ such that $\log a\le G\le-\log b$ pointwise and
$\sup_R(\av{|G-G_R|^2}R)^{1/2}\le C_*(d)$. We use $\log 0=-\infty$.
\end{lemma}
\begin{proof}
We construct positive functions $U,V$ satisfying
$\mathcal L(U,V)\le\Lambda(U,V)$ for a constant $\Lambda$ depending
only on $d$, where
\[
 \mathcal M=\max_{1\le i\le d}M_i,\qquad
 \mathcal L(U,V)=(\mathcal MU+D_aV,\mathcal MV+D_bU).
\]
Here the maximum defining $\mathcal M$ is pointwise. The logarithm of
$U/V$ will have controlled mean oscillation, and projecting a multiple
of it onto the required interval will give $G$.

Each application of $\mathcal L$ is attained by a matrix
$\bigl(\begin{smallmatrix}P&D_a\\D_b&Q\end{smallmatrix}\bigr)$
with policies $P,Q$ chosen on its two inputs. Expand a length-$m$ matrix
product by summing over the choice of component at each step.
The off-diagonal factors alternate. Pair each $D_a$ in written order
with the next $D_b$ if present. Each resulting block $D_a P D_b$
with $k$ diagonal policies uses $k+2$ factors and has norm at most
$C_d^{k+1}$ by Proposition~\ref{weighted}. Each leftover diagonal
policy has uniformly bounded norm. Together these factors contribute
at most $C_d^m$, after increasing $C_d$ if needed. A factor $D_b$ can
remain unpaired at the left end of the off-diagonal sequence, and a
$D_a$ at the right end. Their norms contribute at most
$\Gamma_{\rm grid}=\max(1,\|D_a\|)\max(1,\|D_b\|)$ once per path.
Summing over at most $2^{m+1}$ paths bounds the matrix product norm
by $\Gamma_{\rm grid}C_0(d)^m$ for $m\ge1$.

Each iterate on $\mathbf e=(\one,\one)$ is attained by such a
product. Take $\Lambda=2\max(1,C_0(d))$. The series
$(U,V)=\sum_{m\ge0}\Lambda^{-m}\mathcal L^m\mathbf e$ then converges
on each grid. Positivity, homogeneity, subadditivity and continuity give
$U,V\ge1$ and $\mathcal L(U,V)\le\Lambda(U,V)$. The constant $\Lambda$
depends only on $d$; the grid-dependent factor $\Gamma_{\rm grid}$
enters only the convergence.

Since $M_i U\le\Lambda U$ for every $i$, fixing the side intervals of $R$
and averaging successively bounds $\av UR$ by $\Lambda^d U(z)$ at every
$z\in R$. Applying this also to $V$ gives
$\av UR\le \Lambda^d\inf_R U$ and $\av VR\le\nobreak \Lambda^d\inf_R V$.
For $w=U/V$ these bounds give
\[
 \av wR\av{w^{-1}}R
 \le \frac{\av UR}{\inf_R U}\frac{\av VR}{\inf_R V}
 \le \Lambda^{2d}.
\]
With $u=\log w$, Jensen gives $e^{-u_R}\le\av{w^{-1}}R$ and
$e^{u_R}\le\av wR$. Hence each $\av{e^{\pm(u-u_R)}}R$
is at most $\Lambda^{2d}$. Since $s^2\le e^s+e^{-s}$, the grid
rectangle seminorm of $u/\beta$ is at most
$\sqrt{2}\Lambda^d/\beta$.

The inequalities $a^\beta V\le\Lambda U$ and
$b^\beta U\le\Lambda V$ place $u/\beta$ within $(\log\Lambda)/\beta$
of the nonempty interval $[\log a,-\log b]$ at each cell. Define $G$ by
projecting $u/\beta$ onto this interval at each cell, with an infinite
endpoint wherever its weight vanishes. Thus
$\|G-u/\beta\|_\infty\le(\log\Lambda)/\beta$.
On each rectangle the standard deviation of $G-u/\beta$
is at most $\|G-u/\beta\|_\infty$. The triangle inequality therefore gives
$\|G\|_{\BMO,\mathrm{grid}}\le
C_*=(\sqrt2 \Lambda^d+\log\Lambda)/\beta$.
\end{proof}

\subsection{A correction on all rectangles}
For real $f$ and admissible $\alpha$, write $K_\alpha=K_\alpha(f)$ and use
the mesh partitions of $Q_n$ from Section~\ref{ordinary}. Let $f_n$ be
the cell averages of $f$. On each aligned rectangle $R$, cellwise Jensen
gives $\av{e^{\pm\alpha f_n}}R\le K_\alpha e^{\pm\alpha f_R}$. Thus the
weights $a=e^{\alpha f_n}/K_\alpha$ and $b=e^{-\alpha f_n}/K_\alpha$
satisfy $\av aR\av bR\le1$. Lemma~\ref{finitebarrier} gives a grid
function $G_n$ with $|G_n-\alpha f_n|\le\log K_\alpha$ and quadratic
mean oscillation at most $C_*$ on grid rectangles. Identify $G_n$ with
the function on $Q_n$ that is constant on each cell. On aligned
rectangles the counting and Lebesgue averages agree, so the same
oscillation bound holds there.

The moment hypothesis gives $f\in L^2_{\mathrm{loc}}$, and dyadic
conditional expectation gives $f_n\to f$ in $L^2$ on each fixed
rectangle with dyadic-rational endpoints. Extend the corrections
$h_n=\alpha f_n-G_n$ by zero off $Q_n$. They are bounded by
$\log K_\alpha$, so by separability of $L^1$ a single subsequence
converges weak-$*$ in $L^\infty(\mathbb R^d)$ to a function $h_*$ with
the same bound. On each rectangle $R$ with dyadic-rational endpoints,
$L^2(R)$ tests extend by zero to $L^1(\mathbb R^d)$. Since $R$ is
eventually aligned, $G_n=\alpha f_n-h_n$ converges weakly there in $L^2$
to the residual $\alpha f-h_*$ along the same subsequence. Subtract the
mean and use weak lower semicontinuity to get the oscillation bound
$C_*$ on $R$.

For an arbitrary bounded rectangle, approximate its endpoints by
dyadic rationals in a fixed compact region. Absolute continuity of
the integrals of the residual and its square gives convergence of the
first and second moments. The variance formula transfers the bound
$C_*$ to this rectangle. Thus $h=h_*/\alpha$ satisfies
\[
 \|h\|_\infty\le \alpha^{-1}\log K_\alpha(f),\qquad
 \|f-h\|_{\BMO}\le C_*/\alpha.
\]
By Cauchy--Schwarz this gives the real-valued correction estimate in
Theorem~\ref{distance} with $C_d=C_*$.

\section{Exponential integrability and the converse}\label{converse}
We first obtain exponential integrability on an interval, then
iterate it across coordinates. Applying the resulting rectangle
estimate to $f-h$ gives the lower distance bound.

\paragraph{Intervals.}
Use the stopping argument of John--Nirenberg \cite{jn}. Suppose $r$
is real and has mean absolute oscillation at most $\omega>0$ on
intervals. Bisect a fixed interval $I$ repeatedly, stopping at maximal
dyadic subintervals $J$ with $\av{|r-r_I|}J>2\omega$. Their total length
is at most $|I|/2$. The dyadic parent of each $J$ has average at most
$2\omega$, so $|r_J-r_I|\le4\omega$. Off their union $|r-r_I|\le2\omega$
almost everywhere by differentiation. Repeat in each $J$ with its own
mean. The union of generation $m$ has length at most $2^{-m}|I|$. On its
difference with the next generation, $|r-r_I|\le(4m+2)\omega$. The
intersection of all generations is null. For $c=(\log 2)/8$ it follows
that
\[
 \av{e^{c|r-r_I|/\omega}}I \le e^{2c}\sum_{m\ge0}(e^{4c}/2)^m=:B_0<\infty .
\]

\paragraph{Rectangles.}
The interval estimate is the base case for induction in $d$. Suppose
$r$ is real and has mean absolute oscillation at most $\omega>0$ over
rectangles. For almost every $y\in\mathbb R^{d-1}$ the slice $r(\cdot,y)$
has interval BMO norm at most $\omega$. Indeed, for a rational-ended
interval $I$ write $m_I(y)=\av{r(\cdot,y)}I$. The bound
$\av{|r-q|}{I\times Q}\le \omega+|(m_I)_Q-q|$
holds for every rational $q$ and every transverse cube $Q$. Shrinking $Q$
to $y$, differentiation and countability give, on one full-measure set,
$\av{|r(\cdot,y)-q|}I\le \omega+|m_I(y)-q|$ for all such $I,q$.
Taking $q\to m_I(y)$ and approximating interval endpoints gives the
slice bound on every interval. For fixed $I$, the mean $m_I$ has
$(d-1)$-dimensional rectangle oscillation at most $\omega$ by Jensen's
inequality on each transverse rectangle $Q$,
\[
 \av{|m_I-(m_I)_Q|}Q\le\av{|r-r_{I\times Q}|}{I\times Q}\le \omega.
\]

On $I\times Q$, the difference $r-r_{I\times Q}$ is the sum of
$r-m_I$ and $m_I-(m_I)_Q$. For almost every $y\in Q$, the interval
estimate gives
\[
 \frac{1}{|I|}\int_I e^{c|r(x,y)-r_{I\times Q}|/\omega}\,dx
 \le B_0 e^{c|m_I(y)-(m_I)_Q|/\omega}.
\]
Averaging over $Q$ and using the induction hypothesis bounds the
rectangle moment by $B_0^d$ at exponent $c/\omega$.

\paragraph{Distance to bounded functions.}
For real $f$ and a bounded real $h$, let $r=f-h$ and
$\omega=\|r\|_{\mathrm{BMO}_{1,\mathcal R}}$.
If $\omega>0$, adding $h-h_R$ to $r-r_R$ increases the moment at
exponent $c/\omega$ by at most the factor $\exp(2c\|h\|_\infty/\omega)$.
Hence $a_*(f)\ge c/\omega$. If $\omega=0$, the residual is constant almost
everywhere by testing nested rectangles, so $a_*(f)=\infty$.
In particular, taking $h=0$ shows $a_*(f)>0$ for every real
$f\in\mathrm{BMO}_{1,\mathcal R}$.

Taking the infimum over $h$ gives the lower distance bound with constant
$c$. For the upper bound, take the infimum of $C_*/\alpha$ over
admissible $\alpha$ in the correction estimate. The two bounds give
distance zero exactly when $a_*(f)=\infty$, that is, by downward
closure, exactly when every $\alpha>0$ is admissible.

\paragraph{Complex-valued functions.}
Write $f=u+iv$. Since the modulus of each component of $f-f_R$ is
bounded by $|f-f_R|$, we have $K_\alpha(u),K_\alpha(v)\le K_\alpha(f)$.
Apply the real correction estimate with constant $C_*$ from
Section~\ref{barrier} to obtain $h_u,h_v$, and set $h=h_u+ih_v$.
Squaring the modulus adds the squares of the two components, both for
$h$ and for the centered residual on each rectangle. Thus
\[
 \|h\|_\infty\le \frac{\sqrt{2}}{\alpha}\log K_\alpha(f),\qquad
 \|f-h\|_{\BMO}\le \frac{\sqrt{2}C_*}{\alpha}.
\]
Cauchy--Schwarz gives the same bound for the mean absolute oscillation.

For the lower bound let $h$ be any bounded complex-valued function,
$r=f-h$, and $\omega=\|r\|_{\mathrm{BMO}_{1,\mathcal R}}$. If $\omega>0$,
both real components of $r$ have mean absolute oscillation at most
$\omega$. Using
$|r-r_R|\le |\operatorname{Re}(r-r_R)|+|\operatorname{Im}(r-r_R)|$
and the Cauchy--Schwarz inequality gives
\[
 K_{c/(2\omega)}(r)
 \le \bigl(K_{c/\omega}(\operatorname{Re}r)
           K_{c/\omega}(\operatorname{Im}r)\bigr)^{1/2}
 \le B_0^d.
\]
Adding back $h$ preserves admissibility of $c/(2\omega)$. If $\omega=0$,
both components of $r$ are constant almost everywhere and
$a_*(f)=\infty$. Taking $h=0$ gives $a_*(f)>0$ for complex BMO
functions. Taking the infimum over $h$ gives the lower distance bound.
Taking the infimum over admissible $\alpha$ in the correction estimate
gives the upper bound. Thus
$c_d=c/2$ and $C_d=C_*$ work in Theorem~\ref{distance} for both real and
complex $f$; the real distance bound even holds with $C_d$ in place of
$\sqrt2C_d$. The two distance bounds and downward closure of the
admissible exponents give the zero-distance criterion.

\section*{Acknowledgment}
ChatGPT 6.2 Astra was used under the authors' direction in developing the
construction. AI tools also assisted with proof review and editing.

\end{document}